\documentclass[letterpaper, 10 pt, conference]{ieeeconf}  
\IEEEoverridecommandlockouts                              
\newcommand{\subparagraph}{}
\usepackage{titlesec}

\let\labelindent\relax
\IEEEoverridecommandlockouts
\usepackage{gametix}
\usepackage{mathematix}
\usepackage{mytodos}

\newif\ifArxiv
\Arxivtrue

\acrodef{CM}{center of mass}
\Crefname{figure}{Fig.}{Figs.}
	
\begin{document}
\title{
    Learning MPC for Interaction-Aware Autonomous Driving: 
	A Game-Theoretic Approach
        \thanks{
            This work was supported by the Research Foundation Flanders (FWO) research projects G0A0920N, G086518N, G086318N, and PhD grant 1183822N;
            Ford KU Leuven Research Alliance Project KUL0075;
            Research Council KU Leuven C1 project No. C14/18/068;
            and the Fonds de la Recherche Scientifique -- FNRS and the Fonds Wetenschappelijk Onderzoek -- Vlaanderen under EOS project no 30468160 (SeLMA).
        }
}

\author{%
	\texorpdfstring{%
		Brecht Evens\thanks{
			KU Leuven, Department of Electrical Engineering ESAT-STADIUS -- %
			Kasteelpark Arenberg 10, bus 2446, B-3001 Leuven, Belgium
			\newline
			{\sf
				\{%
					\href{mailto:brecht.evens@kuleuven.be}{brecht.evens},%
					\href{mailto:mathijs.schuurmans@kuleuven.be}{mathijs.schuurmans},%
					\href{mailto:panos.patrinos@kuleuven.be}{panos.patrinos}%
				\}%
				\href{mailto:brecht.evens@kuleuven.be,mathijs.schuurmans@kuleuven.be,panos.patrinos@kuleuven.be}{@kuleuven.be}%
			}%
		}
		\and
        Mathijs Schuurmans
		\and
		Panagiotis Patrinos
	}{%
        Brecht Evens, Mathijs Schuurmans, and Panagiotis Patrinos
	}%
}

\maketitle

\begin{abstract}
We consider the problem of interaction-aware motion planning for automated vehicles in general traffic situations. 
We model the interaction between the controlled vehicle and surrounding road users using a generalized potential game, in which 
each road user is assumed to minimize a common cost function subject 
to shared (collision avoidance) constraints. 
We propose a quadratic penalty method to deal with the shared constraints 
and solve the resulting optimal control problem online using an Augmented Lagrangian method based on PANOC. 
Secondly, we present a simple methodology for learning preferences and constraints of other road users online, based on observed behavior.
Through extensive simulations in a highway merging scenario, we demonstrate the practical efficacy of the overall approach as well as the benefits of the proposed online learning scheme.


\end{abstract}

\section{Introduction}
\label{sec:introduction}
\subsection{Background and motivation}

Despite many efforts in recent years, safe autonomous navigation in the presence of humans remains a highly challenging task.
Traditionally, this problem is tackled by separating 
it into a forecasting problem and a 
planning problem.
That is, predictions of surrounding vehicles' trajectories are generated based on assumptions about their
drivers, and these predictions are then fixed and fed into a planning module.
However, since the mutual interactions between the autonomous vehicle and the surrounding traffic are neglected, this type of formulation may lead to overly conservative behavior,
also referred to as the \emph{frozen robot problem} \cite{Trautman2010}. 
In dense traffic, for instance, such a vehicle will typically be unable to merge into lanes cross condensed junctions, even though in practice, feasible solutions may very well exist.

Early work on social navigation such as \cite{Trautman2010} proposed to solve the frozen robot problem by modelling vehicles as a team of cooperative players,
working together to ensure that each player reaches their destination while avoiding collisions. Autonomous navigation can be related to other problems in game theory as well, such as antagonistic and Stackelberg (i.e., leader-follower) games. Antagonistic games are commonly used to account for the worst-case behavior of the surrounding traffic, resulting in safe, robust control strategies \cite{lachner2000}, which may again suffer from the frozen robot problem. 
Stackelberg games can be used to model interactions where one driver, i.e., the leader, dominates the decision process \cite{Sadigh2018,Liniger2019}.
Despite its computational benefits, a drawback of this approach is that it places the burden of collision avoidance unilaterally on the followers.
This may lead to undesirably aggressive behavior, as the leader (typically the autonomous vehicle) will assume that the followers can anticipate its planned trajectory accurately.


In this work, we aim to account for the mutual interactions between the 
controlled vehicle and other road users by modeling it as a cooperative game with coupled constraints,
following the previously mentioned work on social navigation.
The problem of finding a solution, i.e., a Nash equilibrium, of a game with coupled constraints is called a \ac{GNEP} \cite{Debreu1952} and is in general notoriously challenging to solve. 
Efficiently computing solutions to \acp{GNEP} in real-time poses one of the major challenges of this paper.

Moreover, a game-theoretic formulation of road user interaction requires the controlled vehicle to have access to the (implicit) objective function and constraints of the other road users. Among autonomous vehicles, such knowledge can reasonably be assumed to be obtained via vehicle-to-vehicle communication. However, when human drivers are considered, it needs to be inferred from observations, taking into account the possibly time-varying nature of their actions.
In the unconstrained setting, this is commonly performed by parametrizing the human objective function and then using Inverse Reinforcement Learning to update these parameters \cite{Levine2012,Kuderer2015,Sadigh2018}.
In this work, we aim to design a practical methodology able to also incorporate and learn constraints.

Our contributions are twofold:
\paragraph{A penalty method for interaction-aware motion planning}
	We formalize driving interaction with mutual collision avoidance as \acp{GPG} \cite{Monderer1996}, a subclass of \acp{GNEP}, for which we can find a normalized Nash equilibrium \cite{Rosen1965} by solving a single nonlinear \ac{OCP} directly, instead of having to use dedicated solvers such as \cite{fridovich-keil2020,cleach2020} for general \acp{GNEP}.
	This kind of Nash equilibrium is particularly suited for socially aware motion planning as for such an equilibrium the Lagrange multipliers associated with the shared constraints are equal among all players, which introduces some notion of fairness among the players. This idea has been used in various recent papers on traffic control, i.e., without human drivers, such as \cite{Dreves2018,britzelmeier2021}, where ellipsoidal constraints are used to enforce safety.
	As ellipsoidal overapproximations seem less suitable to model human behavior and may furthermore be conservative in small distances (where interaction is more pronounced), we instead extend the collision avoidance formulation of \cite{Hermans2018}, which involves no approximation of rectangular obstacles.

\paragraph{Online learning based on an optimal control model}
	We propose an online learning methodology to continually update our estimate of other player's costs and constraints by adapting standard inverse optimal control methodologies such as \cite{Hatz2012, menner2020, peters2021} to our game-theoretical framework.
	We show empirically that this approach performs well in practice, demonstrating successful closed-loop navigation, where the certainty-equivalent controller exhibits suboptimal or even dangerous behavior.

\subsection{Notation}
Given two nonnegative integers $a \le b$, let $\N_{[a,b]} \dfn \{n \in \N\; |\; a \le n \le b\}$.
We define $\plus{\wc} \dfn \max\{0,\wc\}$, where $\max$ is interpreted element-wise.
Given a set $U$, we denote $U^N \dfn U \times \dots \times U$ as the $N$-fold Cartesian product with itself, for $N \in \N \setminus \{0\}$.
Given a positive definite matrix $Q \succ 0$, 
let $\|x\|_{Q}^2 = \trans{x}Q x$ denote the weighted square norm.

\section{Problem statement}\label{sec:problem_statement}
We consider the task of controlling an autonomous vehicle interacting with human drivers (and optionally also other autonomous vehicles) in general traffic environments,
and formulate this task as a repeated game embedded in a model predictive control (MPC) scheme.

Formally, a game consists of $M$ players, labeled by the index set $\agents \dfn \{ 1,2,\dots,M \}$. Following the notation in \cite{Sadigh2018}, we denote the index set of autonomous vehicles using symbol $\robot \subseteq \agents$, and of human drivers using symbol $\human = \agents \setminus \robot$.
Each player $\nu \in \agents$ decides on their control variables
$\ups[k]{\nu} \in \Re^{n_\nu}$ at time step $k$.
For notational convenience, we also introduce $\ups[k]{-\nu} \dfn (\ups[k]{i})_{i \in \agents \setminus \{\nu\}}$ as the $(M-1)$-tuple formed by the input variables of all players except for player $\nu$ at time step $k$ and $\us[k] \dfn (\ups[k]{i})_{i\in\agents}$ as the $M$-tuple formed by the input variables of all players at time step $k$.
The state vector $\xps[k]{\nu}$ represents the physical state of player $\nu$, which evolves over time according to the dynamics
$\xps[k+1]{\nu} = \dynp{\nu}\left( \xps[k]{\nu}, \ups[k]{\nu}\right)$.
Finally, $\xs[k]$ denotes the full physical state vectors, i.e., the $M$-tuple of state vectors of all players at time step $k$. 

In order to formulate a game-theoretic description of the control 
task for the autonomous vehicle, we make the following key assumption.
\begin{assumption} \label{assum:driver-cost}
	A human driver $\nu \in \human$ behaves optimally according to a (potentially unknown) receding horizon optimal control cost characterized as a sum of stage costs $\Js{\nu}(\xs[k], \ups[k]{\nu}, \ups[k]{-\nu})$, 
	over a prediction horizon of $N$ time steps.
\end{assumption}
In particular, we assume that at any time step, each player aims to minimize an additive cost 
\(
	\sum_{k=0}^{N-1} \Js{\nu}(\xs[k], \ups[k]{\nu}, \ups[k]{-\nu}) + \Jt{\nu}(\xs[N])
\) 
over an $N$-step prediction horizon, where the stage costs $\Js{\nu}$ may be functions of the state and control variables of all other players. 
Furthermore, we impose hard constraints to model physical (actuation) limits and to ensure safety. Here, we make the distinction between three types of constraints.
First, input constraints $\ups{\nu} \in \boundsps{\nu}$ representing the physical limitations of the vehicle actuation. For convenience, we will assume that $\boundsps{\nu}$ represents box constraints, i.e., 
\begin{equation}\label{eq:box-constraints-step}
	\boundsps{\nu} = \left\{ u \in \Re^{n_\nu} \mid \lbx{b}^{\nu} \leq u \leq \ubx{b}^{\nu} \right\}.
\end{equation}
Second, continuously differentiable player-specific constraints $\hs{\nu}(\xps[k]{\nu}, \ups[k]{\nu}) \leq 0$, such as constraints ensuring that the vehicle remains within the boundaries of the road section, and third, collision avoidance constraints $\Hs(\xs[k], \us[k]) \leq 0$, shared by all players.

Since the players are unaware of the cost function and the constraints of the other players, the formulation up to this point corresponds to the formulation of a partially observable repeated game.
In order to obtain a more tractable formulation, we introduce the following two additional assumptions.
\begin{assumption} \label{assum:knowledge}
	The autonomous vehicles can either communicate with each other or are in possession of the cost function and the constraints of each other.
	\label{ass:robot_cost_function_known}
\end{assumption}
\begin{assumption} \label{assum:parametric}
	The cost function and constraints of each human driver $\nu \in \human$ can be parametrized by
	$\params{\nu}$, and the resulting parametrized cost function and constraints are representative of the real behavior of this human driver.
	\label{ass:human_cost_function_known}
\end{assumption}
Note that \Cref{assum:knowledge} could be satisfied in practice rather easily.
More restrictive is \Cref{assum:parametric}, first proposed by Kuderer et al. \cite{Kuderer2015}, as it assumes that the full range of possible driving styles can be obtained by varying the parameters $\params{\nu}$. Therefore, we need to use a sufficiently expressive function class to represent the human's cost and constraints, such that in principle arbitrarily realistic behavior can be obtained.
Our approach towards this parametrization is detailed in \Cref{sec:online_learning}.
Under previous assumptions, the game reduces to a fully observable repeated game where each player $\nu \in \agents$ is assumed to be solving a coupled discrete-time finite-horizon optimal control problem, described by
\begin{mini}[2]{\substack{\xps[\secondindex]{\nu},\hdots,\xps[\endindex]{\nu},\\\ups[\initindex]{\nu},\hdots,\ups[\endindex-1]{\nu}\in \boundsps{\nu}}}
	{\hspace{-6pt}\sum_{k=\initindex}^{\endindex-1} \Js[k]{\nu}(\xs[k], \ups[k]{\nu}, \ups[k]{-\nu}) + \Jt{\nu}(\xs[\endindex])}
	{\label{mini:discrete_time_optimal_control_problem_multiple_players}}
	{\hspace{-5pt}\problem{\nu}:\hspace{-11pt}}
	\addConstraint{\hspace{-7pt}\xps[k+1]{\nu} = \dynp[k]{\nu}(\xps[k]{\nu}, \ups[k]{\nu}), \hspace{6pt}}{}{\forall k \in \N_{[\initindex,\endindex-1]},}
	\addConstraint{\hspace{-7pt}\hs{\nu}(\xps[k]{\nu}, \ups[k]{\nu}) \leq 0, \hspace{6pt}}{}{\forall k \in \N_{[\initindex,\endindex-1]},}
	\addConstraint{\hspace{-7pt}\hst{\nu}(\xps[N]{\nu}) \leq 0,}
	\addConstraint{\hspace{-7pt}\Hs(\xs[k], \ups[k]{\nu}, \ups[k]{-\nu}) \leq 0, \hspace{6pt}}{}{\forall k \in \N_{[\initindex,\endindex-1]},}
	\addConstraint{\hspace{-7pt}\Hst(\xs[\endindex]) \leq 0.}
\end{mini}
Such a system of coupled optimization problems $\{\problem{\nu}\}_{\nu\in \agents}$ is called a \acs{GNEP}, which as mentioned before is in general very challenging to solve due to the coupling through both the cost function and the constraints.
In the next section, we will impose some additional structure onto the given problem, which will facilitate an efficient solution procedure.

\section{Solution of the game formulation}
\subsection{Reformulation to standard game-theoretic form}
In order to simplify notation in the sequel, we 
reformulate \eqref{mini:discrete_time_optimal_control_problem_multiple_players} into a more compact form. To this end, 
we eliminate the state variables
by defining $\Dynp[k]{\nu}$, $k \in \N$ as the solution maps to the dynamics $\dynp{\nu}$ starting from $x_0$, i.e., $\Dynp[0]{\nu}\left(\upf{\nu}\right) \dfn \xps[0]{\nu}$ and
\begin{align*}
	\Dynp[k+1]{\nu}\left(\upf{\nu}\right) &\dfn \dynp{\nu}\left(\Dynp[k]{\nu}\left(\upf{\nu}\right), \ups[k]{\nu}\right),\\
	\Dyn{k}\left(\upf{}\right) &\dfn \left(\Dynp[k]{1}\left(\upf{1}\right),\dots,\Dynp[k]{M}\left( \upf{M}\right)\right),
\end{align*}
where $\upf{\nu} \dfn \left(\ups[\initindex]{\nu},\dots, \ups[\endindex-1]{\nu}\right)$ and $\uf \dfn \left(\us[\initindex],\dots, \us[\endindex-1]\right)$.
We may then define the corresponding single shooting cost function $\J{\nu}$, the player-specific constraints $\hf{\nu} \equiv (\hs{\nu}_{k})_{k=\initindex}^{\endindex}$ and shared constraints 
$\Hf \equiv (\Hs_{k})_{k=\initindex}^{\endindex}$ as
\begin{align*}
	\J{\nu}\left( \upf{} \right) 
		&\dfn 
			\textstyle\sum_{k=\initindex}^{\endindex-1} \Js{\nu}\left(\Dyn{k}( \upf{}\right), \ups[k]{}) 
			+ \Jt{\nu}\left(\Dyn{N}\left(\upf{\nu} \right)\right),
	\\
	\hs{\nu}_{k}\left(\upf{\nu}\right) 
		&\dfn \hs{\nu}\left(\Dynp[k]{\nu}\left(\upf{\nu}\right), \ups[k]{\nu}\right),\;
	\hs{\nu}_{\endindex} \left(\upf{\nu}\right) \dfn \hst{\nu}\left(\Dynp[k]{\nu}\left(\upf{\nu}\right)\right), \\ 
	\Hs_k(\upf{})
		&\dfn \Hs \left(\Dyn{k}\left(\upf{}\right), \ups[k]{}\right), \;  
	\Hs_{\endindex}(\upf{}) \dfn \Hst \left(\Dyn{N}\left(\upf{}\right)\right),
\end{align*}
for $k \in \N_{[\initindex, \endindex]}$. We furthermore introduce the constraint set
\begin{equation}
	\Uf^\nu(\upf{-\nu}) = 
	\left\{\upf{\nu} \in \boundspf{\nu}\middle\vert\ \Hf(\uf) \leq 0, \ \hf{\nu}(\upf{\nu}) \leq 0
	\right\},
	\label{eq:GNEP_constraint_set}
\end{equation}

where $\upf{-\nu} \dfn \left(\ups[\initindex]{-\nu},\dots, \ups[\endindex-1]{-\nu}\right)$ and $\boundspf{\nu} \dfn (\boundsps{\nu})^N$.
Thus, \eqref{mini:discrete_time_optimal_control_problem_multiple_players}
may be compactly written in the 
standard GNEP form
\begin{mini}[2]{\upf{\nu} \in \Upf{\nu}(\upf{-\nu})}
	{\J{\nu}(\upf{\nu}, \upf{-\nu}), \quad \nu \in \agents. }
	{}{}
	\label{mini:GNEP}
\end{mini}
We now formally define the sought solution.
\begin{definition}[Generalized Nash equilibrium]
	We say that an $M$-tuple $\uf-$ is a \ac{GNE} of \ac{GNEP} \eqref{mini:GNEP} if no player $\nu$ can decrease their cost by unilateral deviation. That is, for all $\nu \in \agents$, 
\begin{equation}
	\J{\nu}(\upf-{\nu}, \upf-{-\nu}) \le \J{\nu}(\upf{\nu}, \upf-{-\nu}),\quad \forall \upf{\nu} \in \Upf{\nu}(\upf-{-\nu}).
	\label{eq:Nash equilibrium_GNEP}
\end{equation}
\end{definition}

\subsection{Generalized potential games}\label{subsec:Gauss-Seidel}
In the context of traffic interaction, it is 
natural to impose some additional structure on the cost functions of the different interacting players. 
This will allow us to specialize problem \eqref{mini:GNEP} to the form of \iac{GPG} \cite{Monderer1996} --- a subclass of \acp{GNEP} with shared constraints.
\begin{definition}[Generalized potential game] \label{def:GPG}
	A \ac{GNEP} corresponds to a \ac{GPG} if:
	\begin{enumdefinition}
		\item \label{cond: condition 1 GPG}there exists a nonempty, closed set $\Uf \subseteq \Re^n$, with $n \dfn \sum_{\nu = 1}^{M}n_\nu$, such that for every $\nu \in \agents$,
		\begin{equation}
			\Upf{\nu}\left(\upf{-\nu}\right) = \left\{\upf{\nu} \in \boundspf{\nu}\ \mid \left(\upf{\nu}, \upf{-\nu}\right) \in \Uf\right\};
			\label{eq:set-valued mapping}
		\end{equation}
		\item \label{cond: condition 2 GPG}there exists a continuous function $\Pot : \Re^n \rightarrow \Re$ such that for every $\nu \in \agents$ and  for every $\upf{-\nu} \in \prod_{i \in \agents \setminus \{\nu\}} \boundspf{i}$,
		$$
		\tfrac{\partial P(\upf{\nu}, \upf{-\nu})}{\partial \upf{\nu}} = \tfrac{\partial \J{\nu}(\upf{\nu}, \upf{-\nu})}{\partial \upf{\nu}}.
		$$
	\end{enumdefinition}
\end{definition}
\Cref{cond: condition 1 GPG} readily follows from the problem setup. Using \eqref{eq:GNEP_constraint_set}, we define
\ifArxiv
\begin{equation*}
	\Uf = 
	\left\{\uf \in \Re^n\middle\vert\ \begin{aligned}\Hf(\uf) &\leq 0, \\ \hf{}(\uf) &\leq 0\end{aligned}
	\right\},
\end{equation*}
\fi

with $\hf{}(\uf) \dfn \left( \hf{1}(\upf{1}), \dots, \hf{M}(\upf{M}) \right)$. Closedness of $\Uf$ directly follows from continuity of $\hf{}$ and $\Hf$.

\Cref{cond: condition 2 GPG} states that in fact, all players are minimizing the same function $\Pot(\uf)$,
also called the \textit{potential function} of the game.
This trivially holds when 
the cost function of each player $\nu$ can be decomposed as $\J{\nu}(\uf) = c(\uf) + d^\nu(\upf{\nu})$,
where $c(\uf)$ is a common term and $d^\nu(\upf{\nu})$ is a player-specific term, such that $\Pot(\uf) = c(\uf) + \sum_{\nu = 1}^{M}d^\nu(\upf{\nu})$.
In our particular application, we can thus design the cost function of each vehicle as the sum of a player-specific term (e.g., tracking errors, control effort, slew rate, \dots) and a common term (e.g., proximity to other vehicles). Such costs are often used in literature, and they are able to capture most of the behavior of human drivers \cite{Kuderer2015}.

Under these conditions, finding a \ac{GNE} of 
\eqref{mini:discrete_time_optimal_control_problem_multiple_players}
corresponds to solving a single optimal control problem, given by
\ifArxiv
\begin{mini}{\uf \in \boundsf}
	{\Pot(\uf)}{}{\label{mini:gpg}}
	\addConstraint{\Hf(\uf) \leq 0,}{}{}
	\addConstraint{\hf{}(\uf) \leq 0,}{}{}
\end{mini}
\fi 
where $\boundsf \dfn \prod_{i \in \agents} \boundspf{i}$ is the Cartesian product of the player-specific actuator bounds. 

\subsection{NMPC formulation}
\subsubsection{Dynamics}
We will represent the vehicle dynamics using a simple \textit{kinematic bicycle model}\ifArxiv (cfr. \cref{fig:kinematic bicycle model})\fi, i.e.,
\begin{equation}
	\begin{gathered}
		\begin{aligned}
			\dot\px &= v\cos \left( {\psi + \beta } \right) \qquad \dot \psi = \frac{{v}} {{\lr}}{\sin \left( {\beta } \right)}\\
			\dot\py &= v\sin \left( {\psi + \beta } \right) \qquad \dot v = a - \mu v\\
		\end{aligned}\\
	\end{gathered}
	\label{eq:kinematic_bicycle_model}
\end{equation}
with $\beta = \tan ^{ - 1} \left( {\nicefrac{{l_r}} {{(l_f + l_r)}}\tan \delta } \right)$, 
where $\px$ and $\py$ denote the position of the \ac{CM}; $\psi$ the heading angle, $\delta$ the steering angle, $v$ and $a$ respectively the velocity and the acceleration at the \ac{CM}, $l_r$ and $l_f$ the distance between the \ac{CM} and respectively the front and the rear axis, $\mu$ the friction coefficient and $\beta$ the slip angle. The control actions of each driver are the acceleration $a$ and the steering angle $\delta$.
\ifArxiv
\begin{figure}
	\centering
	\includegraphics[width=0.5\columnwidth]{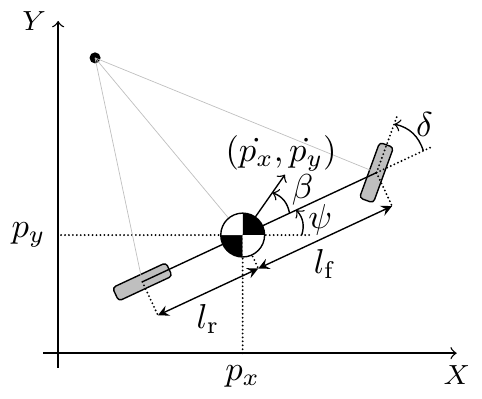}
	\caption{The kinematic bicycle model of \protect \eqref{eq:kinematic_bicycle_model}.}
	\vspace{-0.5cm}
	\label{fig:kinematic bicycle model}
\end{figure}
\fi 
We will restrict our attention to the nominal steering region where $\delta \in (-\frac{\pi}{2}, \frac{\pi}{2})$ for which the dynamics is smooth.

\subsubsection{Collision avoidance constraints}
We model each vehicle as a rectangle, represented as
\begin{equation}
	O^\nu \dfn \left\{ p \in \Re^{2} \mid \trans{p}a_i^\nu  + b_i^\nu > 0, i \in \N_{[1,4]}\right\}.
\end{equation}
Thus, for an arbitrary point $p$,
\begin{equation}
	p \not\in O^\nu \iff \Psi^\nu(p) = \prod_{i=1}^4 \plus{\trans{p} a_i^\nu + b_i^\nu} = 0.
	\label{eq:product_equality_constraint}
\end{equation}
This is a special case of the obstacle avoidance formulation introduced in \cite{Sathya2018,Hermans2018}.
Note that the mapping $\Psi^\nu : \Re^2 \rightarrow \Re$ is such that $\frac{1}{2}\nrm{\Psi^\nu(\cdot)}^2$ is continuously differentiable with Lipschitz-continuous gradient $\nabla_p \frac{1}{2}\nrm{\Psi^\nu(p)}^2$%
\ifArxiv, given by
\begin{equation*}
	\begin{cases}
		\sum_{i=1}^4 (\trans{p} a_i^\nu + b_i^\nu) \prod_{j \neq i} (\trans{p}a_j^\nu + b_j^\nu)^2 a_i^\nu, & \text{if } p \in O^\nu,\\
		0, & \text{otherwise}.
	\end{cases}
\end{equation*}
\fi 
We impose collision avoidance between two vehicles by imposing that each corner of every vehicle lies outside the other vehicle using \eqref{eq:product_equality_constraint}. 
Additionally, we impose that also the nose of each vehicle is also not located within the other vehicle, as this greatly improves performance of our methodology in practice. This yields 10 equality constraints per time step.

\subsubsection{Solver}
We approximately solve \eqref{mini:gpg} by 
treating the shared constraints using the quadratic penalty method with penalty parameters $\Sigma$, such that each subproblem is given by
\ifArxiv
\begin{mini}{\uf \in \boundsf}
	{\Pot(\uf) + \frac{1}{2}\nrm{\Hf(\uf)}_{\Sigma}^2}{}{\label{mini:penalized_gpg}}
	\addConstraint{\hf{}(\uf) \leq 0.}{}{}
\end{mini}
\fi
The player-specific constraints $\hf{}(\uf) \leq 0$ are 
imposed using the \ac{ALM}.
If $\hf{}(\uf)$ is continuously differentiable (as will commonly be the case), then this strategy leads to inner subproblems with a continuously differentiable cost function and simple box constraints, which can be solved approximately using \Panoc{} \cite{Stella2017}.
If this methodology converges to a feasible pair ($\uf-, \bar{y}$), it is an approximate KKT point of \eqref{mini:gpg}.
Otherwise, it converges to an approximate KKT point of the infeasibility \cite{Sopasakis2020}.


\section{Online learning of parameters}\label{sec:online_learning}
We have thus far assumed full access 
to the human cost $\J{\nu}$ and private constraint functions $\hf{\nu}$, $\nu \in \human$.
We now relax this assumption by introducing a practical methodology for learning this information online from observed data.
Our approach is similar in spirit to the methodologies proposed in \cite{Hatz2012,menner2020,peters2021},
although there are some key differences.

First of all, we parametrize not only the cost but also the constraints.
For each human player $\nu \in \human$, we represent the cost and constraint mappings
appearing in \eqref{mini:discrete_time_optimal_control_problem_multiple_players} 
in a parametric form, which we denote
$\Js{\nu}(\wc; \params{\nu}), \Jt{\nu}(\wc; \params{\nu})$,
$\hs{\nu}(\wc; \params{\nu})$,
$\hst{\nu}(\wc; \params{\nu})$
and $\boundsps{\nu}(\params{\nu})$, 
with parameter vector $\params{\nu} \in \Re^{p_\nu}$.
We emphasize that during this design step, care must be taken that \ref{cond: condition 1 GPG} and \ref{cond: condition 2 GPG} hold regardless of the value of $\params{\nu}$.

Second, we assume that the observed control actions $\us^[t]$ of the other players only correspond to the first entries of their respective optimal control sequences $\uf-{}$. This is motivated by the observation that the remaining $N-1$ controls are open-loop predictions which cannot be observed, as they may significantly differ from the closed-loop policy. Considering again the compact representation \eqref{mini:penalized_gpg}, we define the Lagrangian 
\begin{equation}
	\Lag_\Sigma(\uf, y) \dfn \Pot(\uf) + \<y, \mathbf{g}(\uf)\> + \frac{1}{2}\nrm{\Hf(\uf)}_{\Sigma}^2,
\end{equation}
where 
$\mathbf{g}(\uf) \dfn \left(
	\hf{}(\uf),
	\uf{} - \ubx{b}, 
	\lbx{b} - \uf{}
\right)$. Here, we have introduced $\ubx{b}$ and $\lbx{b}$ as
the concatenation of the player-specific actuation bounds $\ubx{b}^\nu, \lbx{b}^\nu$ in \eqref{eq:box-constraints-step}, leading to the 
global box constraints $\boundsf = \{ \uf \in \Re^n \mid \lbx{b} \leq \uf \leq \ubx{b} \}$.
If $\uf-$ is a local minimizer of optimization problem \eqref{mini:penalized_gpg} for given penalty parameters $\Sigma$ and a suitable constraint qualification holds at $\uf-$ then there exists a $\bar{y}$ such that
\ifArxiv
\begin{subequations}
\label{eq:KKT-conditions}
\begin{align}
\nabla_{\uf} \Lag_\Sigma\left(\uf-, \bar{y}; \params{}\right) &= 0,\\
\mathbf{g}(\uf-\-; \params{}) &\leq 0,\\
\bar{y} &\ge 0,\\
\bar{y}\trans{}\mathbf{g}(\uf-; \params{}) &= 0.
\end{align}
\end{subequations}
\fi

The main idea behind the parameter estimation procedure is 
to select parameter values $\params{} \dfn (\params{i})_{i \in \human}$, such that the observed 
behavior approximately matches the optimality conditions \eqref{eq:KKT-conditions}.
This naturally leads to the following parameter estimation procedure: 

Given the current estimate $\params^{}_{t}$ for the parameters, set 
the updated parameters $\params^{}_{t+1}$ as the solution of
\ifArxiv
\begin{mini}[2]{\params{}, \uf{}, y}
	{\left\lVert\nabla_{\uf{}} \Lag_{\widetilde{\Sigma}}\left(\uf, y; \params{}\right)\right\rVert_2^2 + 
	r(\params{}, \params^{}_{t})}
	{\label{mini:human_parameter_estimation}}{}
	\addConstraint{\mathbf{g}(\uf; \params{}) \le 0,}{}
	\addConstraint{y \ge 0,}{}
	\addConstraint{y\trans{}\mathbf{g}(\uf; \params{}) = 0,}{}
	\addConstraint{\us[0]=\us^[t],}{}
\end{mini}
\fi
where 
\(
	r(\params{}, \params^{}_{t}) \dfn \|  \params{} - \params^{}_{t} \|_{\Xi}^2,
\) 
is a regularization term with parameters $\Xi \dfn \diag(\xi) \in \Re^{p\times p}_{\pplus\pplus}$ where $p = \sum_{i \in \human} p_i$, which penalizes large deviations between subsequent parameter estimations. The introduction of this term is motivated by the assumption that the human driving style should only change slowly over time. The penalty parameters $\widetilde{\Sigma}$ are selected equal to the values obtained from the solution of \eqref{mini:gpg} by the controlled vehicle at the previous time step.

Optimization problem \eqref{mini:human_parameter_estimation} updates $\params^{}_t$ based on the information gathered during a single time step $t$ of the game.
Alternatively, an update can be performed based on the previous $L$ observations $\ups^[t]{-\nu},\dots,\ups^[t-L+1]{-\nu}$ by adding the corresponding optimization variables $({\uf}_{(l)}, y_{(l)})$, cost terms $\lVert\nabla_{\uf} \Lag_{\widetilde{\Sigma}_{(l)}}({\uf}_{(l)}, y_{(l)}; \params{})\rVert_2^2$ and constraints to the optimization problem, with ${l = t-L+1, \dots, t}$. 
\ifArxiv
Using only a single observation, i.e., taking $L = 1$, the parameter estimation problem may admit a large number of (approximate) solutions,
as a great deal of freedom remains to choose both the completion of the observed control action over the horizon as well as the parameters $\params{}$.

As $L$ is increased and more state-input pairs (with shared $\params{}$) are taken into account, the ambiguity in the optimal parameter choices is expected to decrease (at the cost of increased computational complexity).
\fi

It is well-known that parameter estimation tasks of this form may inherently suffer from some identifiability issues, as some parameters may be invisible to the learning methodology \cite{soderstrom1983}.
Indeed, there may be a large set of parameter values for which the observed sequence optimizes the \ac{GPG}, locking the learning methodology onto any of these values, regardless of the amount of data \cite[Sec. 6.1.2]{Bertsekas2005}.
Fortunately, our goal is not to arrive at the exact parameters,
but rather at any parameter vector that is consistent with the observed 
driving behavior. The numerical experiments of \Cref{sec:case_studies} show that although parameter estimates in general do not converge to the true values, the closed loop behavior of the controlled vehicle is still significantly improved when using online learning.
\ifArxiv
\begin{remark}
	Optimization problem \eqref{mini:human_parameter_estimation} can be regarded as a hierarchical optimization problem, where the upper level consist of estimating the parameters $\params{}$ and the lower level is an optimal control problem, which can be exploited to design dedicated solution methodologies \cite{Hatz2012}. However, in this work, we opt for the controlled vehicle to solve this optimization problem directly using \Panoc, where we address both the primal feasibility constraints and the complementarity constraints using \ac{ALM}.
	To warm start learning procedure \eqref{mini:human_parameter_estimation} for an autonomous vehicle $\nu$, we use the pair ($\uf-, \bar{y}$) which this vehicle has obtained by solving the \ac{GPG} under its current estimation of the human parameters.
\end{remark}
\fi

\section{Case studies}\label{sec:case_studies}
\subsection{Implementation details}\label{sec:implementation}
All involved optimization problems are implemented using \casadi{} \cite{Andersson2019} in Python and solved in real-time using \panocalm{}\cite{pas2021}
on a 1.7GHz AMD Ryzen 7 PRO 4750U processor with 32 GB RAM (for horizon length $N = 15$ and $\Ts = \SI{0.2}{\second}$).
The source code of the driving simulator and videos of the experiments are available at \href{https://brechtevens.github.io/GPG-drive}{https://brechtevens.github.io/GPG-drive}. 
\subsection{Scenario setup}
\begin{figure}
	\includegraphics[width = \columnwidth]{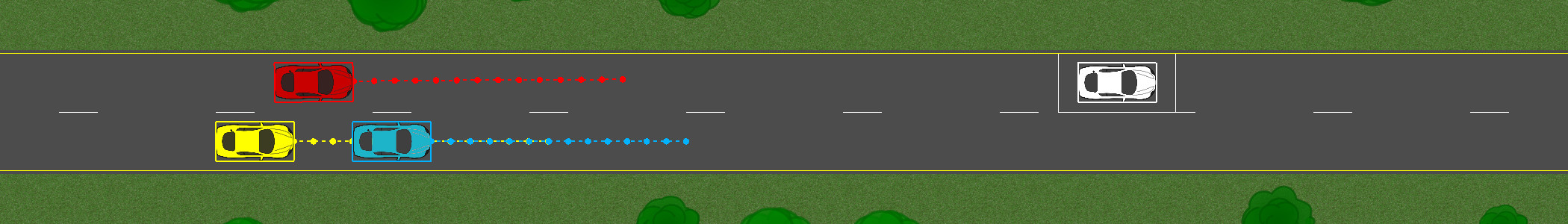}
	\caption{The merging scenario setup.}
	\label{fig:two-dimensional GPG_merging}
	\vspace{-0.5cm}
\end{figure}

We focus on a merging scenario involving three moving vehicles and one stationary vehicle with dimensions $\lr = \SI{2}{\meter}$, $\lf = \SI{2}{\meter}$ and width $\SI{2}{\meter}$.
The initial situation is illustrated in \cref{fig:two-dimensional GPG_merging}, where the controlled vehicle, i.e., the red vehicle $\redplayer$, attempts to merge between the yellow vehicle $\yellowplayer$ and the blue vehicle $\blueplayer$.
We model this scenario as a two-player game as summarized in \Cref{tab:two-dimensional GPG merging dynamics}: the blue vehicle is considered as an obstacle with constant velocity, as is commonly done in traffic modelling on highways \cite{Gipps1986}.
\begin{table}[ht]
	\centering
	\caption{Characteristics and initial states of the vehicles.}
	\label{tab:two-dimensional GPG merging dynamics}
	\begin{tabular}{lcccccc}
		\toprule
		Vehicle & Type & Dynamics & \multicolumn{4}{c}{Initial states}\\
		\cmidrule(lr){4-7}
		& & & $\px[0]$ & $\py[0]$ & $\psi_0$ & $v_0$\\
		\midrule
		Red & \ac{GPG}& Kinematic bicycle & $3$ & $3$ & $0$ & $5$\\
		Yellow & \ac{GPG}& Double integrator & $0$ & $0$ & $0$ & $5$\\
		Blue & Simple & Constant velocity & $7$ & $0$ & $0$ & $5$\\
		White & Simple & Stationary obstacle & $45$ & $3$ & $0$ & $0$\\
		\bottomrule
	\end{tabular}
\end{table}
Cost functions are selected for both players which only contain player-specific cost terms and common cost terms, guaranteeing that \ref{cond: condition 2 GPG} is satisfied. For the controlled vehicle we propose the following player-specific cost terms:
\begin{subequations}
	\begin{align}
	\Js{\redplayer} & \dfn \params{\redplayer}_1(\py[k]^\redplayer)^2 + \params{\redplayer}_2(a_k^\redplayer)^2 + \params{\redplayer}_3({\delta}_k^\redplayer)^2,\\
	\Jt{\redplayer} & \dfn \params{\redplayer}_1(\py[N]^\redplayer)^2.
	\end{align}
\end{subequations}
For the yellow vehicle, we assume a quadratic cost on deviations from the desired distance from the blue vehicle:
\begin{subequations}
	\begin{align}
	\Js{\yellowplayer} & \dfn \params{\yellowplayer}_1(\px[k]^\blueplayer - \px[k]^\yellowplayer - d^\yellowplayer_\mathrm{des})^2 + \params{\yellowplayer}_2(a_k^\redplayer)^2,\\
	\Jt{\yellowplayer} & \dfn \params{\yellowplayer}_1(\px[k]^\blueplayer - \px[k]^\yellowplayer - d^\yellowplayer_\mathrm{des})^2,
	\end{align}
\end{subequations}
where $d^\yellowplayer_\mathrm{des} = 3$.
The ground-truth values for the parameters of the red and yellow vehicle can be found in \Cref{tab:two-dimensional GPG merging scenario cost parameters}, where different parameters are introduced to model either courteous or stubborn behavior of the yellow vehicle.
\begin{table}[ht]
	\centering
	\caption{Cost function parameters.}
	\label{tab:two-dimensional GPG merging scenario cost parameters}
	\begin{tabular}{lcc}
		\toprule
		Vehicle & behavior & Parameters\\
		\midrule
		Red vehicle & - & $\params{\redplayer} = \begin{bmatrix} 0.0.05 & 0.1 & 0.5 \end{bmatrix}\trans{}$\\	
		Yellow vehicle & Stubborn & $\params{\yellowplayer} = \begin{bmatrix} 10 & 0.1 \end{bmatrix}\trans{}$\\
		Yellow vehicle & Courteous & $\params{\yellowplayer} = \begin{bmatrix} 0.02 & 0.1 \end{bmatrix}\trans{}$\\
		\bottomrule
	\end{tabular}
\end{table}
To model the human tendency to keep a certain safety distance, we add Gaussian cost terms penalizing driving near another vehicle \cite{Levine2012, Sadigh2018}, i.e.,
$c(\uf) = \kappa e^{\nicefrac{-d^\top K d}{2}}$, where $d$ denotes the difference between the positions $p$ of two different players, with $\kappa=4$ and $K = \diag(4, 2.25)$.
\ifArxiv
The resulting cost function of the controlled vehicle is visualized in \cref{fig:merging heatmap}.
\begin{figure}
	\includegraphics[width = \columnwidth]{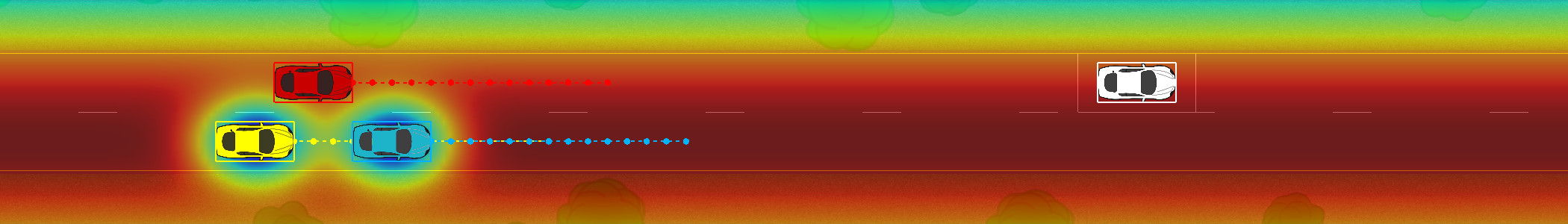}
	\caption{Stage cost for the controlled vehicle.}
	\vspace{-0.25cm}
	\label{fig:merging heatmap}
\end{figure}
\fi 

Both players attempt to avoid collisions and player-specific boundary constraints have been added to the optimal control problem of the red vehicle to ensure that this vehicle remains within the boundaries of the road section. Finally, the inputs of both players are upper and lower bounded. 

\subsection{Numerical results}
\begin{figure}[ht]
	\centering
	\begin{tabular}{wl{0.20\columnwidth} >{\centering\arraybackslash}m{0.69\columnwidth}}
		\toprule
		\hspace{-0.3cm} \scriptsize Belief & \hspace{-2.2cm} \normalsize Courteous vehicle\\
		\midrule
	\end{tabular}
	\begin{subfigure}[t]{0.49\columnwidth}
		\makebox[0.06\columnwidth][r]{\makebox[10pt]{\raisebox{15pt}{\rotatebox[origin=c]{90}{\scriptsize Courteous}}}}%
		\adjincludegraphics[width = \columnwidth-1pt, trim={0 0 {.5\width} 0},clip, cfbox=green 1pt 0pt]{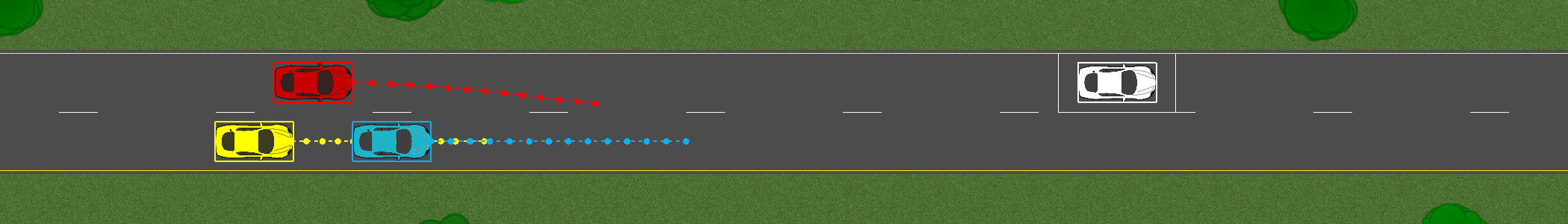}
		\makebox[0.06\columnwidth][r]{\makebox[10pt]{\raisebox{15pt}{\rotatebox[origin=c]{90}{\scriptsize Stubborn}}}}%
		\adjincludegraphics[width = \columnwidth-1pt, trim={0 0 {.5\width} 0},clip,cfbox=orange 1pt 0pt]{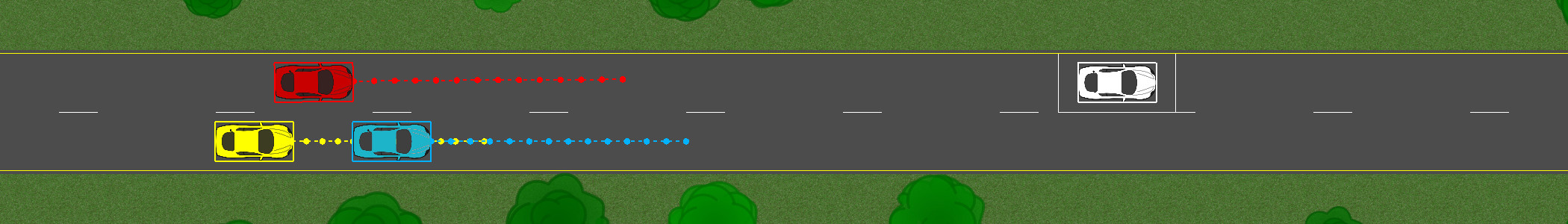}
	\end{subfigure}
	\hspace{-0.015\columnwidth}
	\begin{subfigure}[t]{0.49\columnwidth}
		\adjincludegraphics[width = \columnwidth-1pt, trim={0 0 {.5\width} 0},clip, cfbox=green 1pt 0pt]{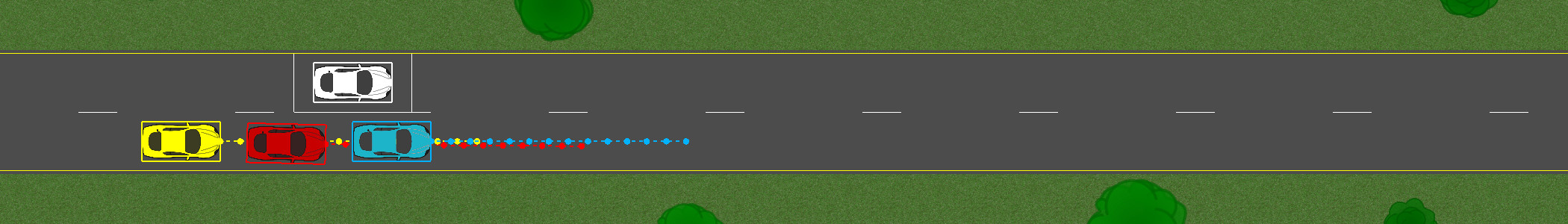}
		\adjincludegraphics[width = \columnwidth-1pt, trim={0 0 {.5\width} 0},clip, cfbox=orange 1pt 0pt]{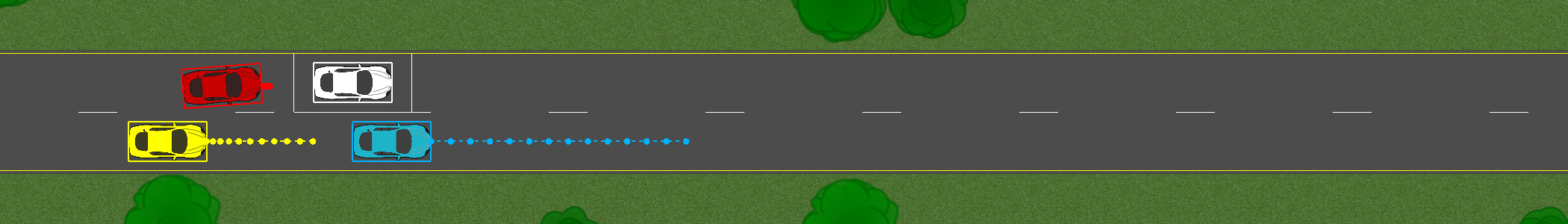}
	\end{subfigure}
	\phantom{\tiny abc}
	\begin{tabular}[ht]{m{0.20\columnwidth} >{\centering\arraybackslash}m{0.69\columnwidth}}
		\toprule
		\hspace{-0.3cm} \scriptsize Belief & \hspace{-2.2cm} \normalsize Stubborn vehicle\\
		\midrule
	\end{tabular}
	\begin{subfigure}[t]{0.49\columnwidth}
		\makebox[0.06\columnwidth][r]{\makebox[10pt]{\raisebox{15pt}{\rotatebox[origin=c]{90}{\scriptsize Stubborn}}}}%
		\adjincludegraphics[width = \columnwidth-1pt, trim={0 0 {.5\width} 0},clip, cfbox=green 1pt 0pt]{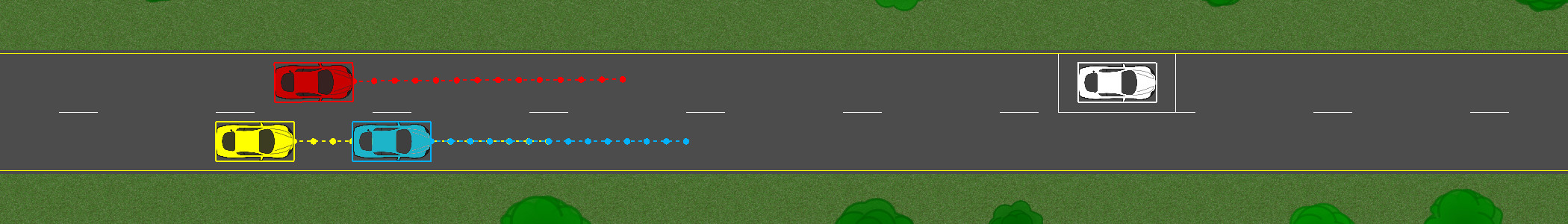}
		\makebox[0.06\columnwidth][r]{\makebox[10pt]{\raisebox{15pt}{\rotatebox[origin=c]{90}{\scriptsize Courteous}}}}%
		\adjincludegraphics[width = \columnwidth-1pt, trim={0 0 {.5\width} 0},clip, cfbox=red 1pt 0pt]{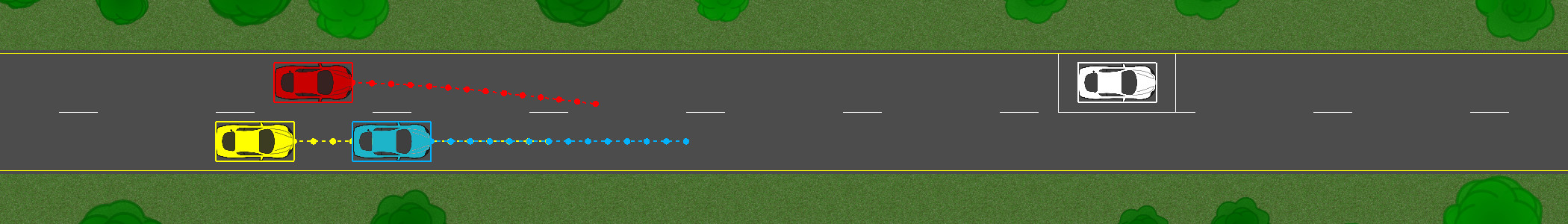}
		\caption{Time step 1}
	\end{subfigure}
	\hspace{-0.015\columnwidth}
	\begin{subfigure}[t]{0.49\columnwidth}
		\adjincludegraphics[width = \columnwidth-1pt, trim={0 0 {.5\width} 0},clip, cfbox=green 1pt 0pt]{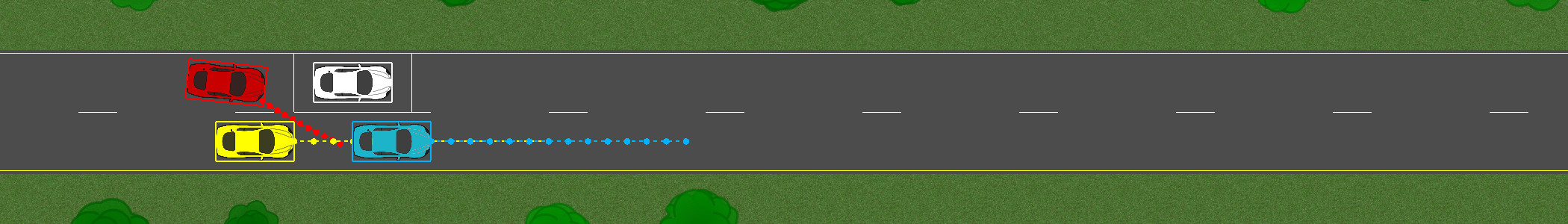}
		\adjincludegraphics[width = \columnwidth-1pt, trim={0 0 {.5\width} 0},clip, cfbox=red 1pt 0pt]{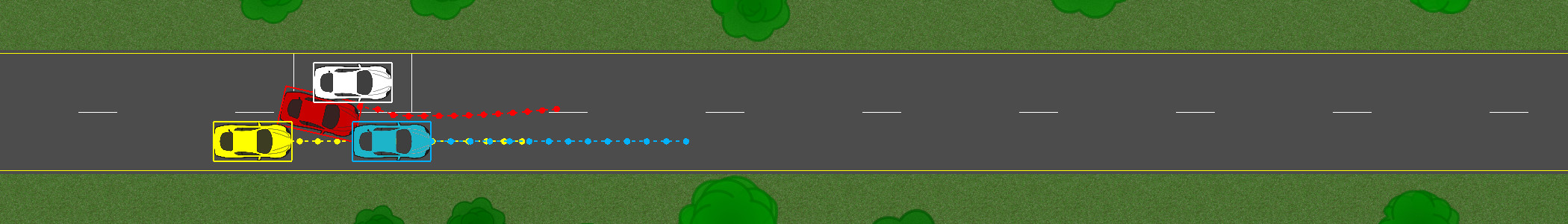}
		\caption{Time step 40}
	\end{subfigure}
	\caption{The planned trajectories at respectively the 1\textsuperscript{st} and 40\textsuperscript{th} time step of the game without online learning.}
	\vspace{-0.4cm}
	\label{fig:two-dimensional GPG merging experimental results}
\end{figure}

The introduced \ac{GPG} formulation is solved for the first $55$ time steps.
We solve \eqref{mini:gpg} for feasibility $10^{-2}$ and tolerance $10^{-3}$.
The yellow vehicle has a perfect estimation of the parameters of the controlled vehicle, but the controlled vehicle might have a wrong belief about the parameters $\params{\yellowplayer}$.
The yellow vehicle either behaves courteously or stubbornly (cf. \Cref{tab:two-dimensional GPG merging scenario cost parameters}).
In the first run, we simply set the parameter estimates of the controlled vehicle fixed to either of the corresponding sets of values (we refer to this as its \emph{belief}).
The obtained closed loop simulations corresponding to these 4 different combinations are visualized in \cref{fig:two-dimensional GPG merging experimental results} at different time steps.
\captionsetup[subfigure]{labelformat=empty}
\captionsetup[subfigure]{labelformat=default}
Consider the scenario where the yellow vehicle behaves courteously (\cref{fig:two-dimensional GPG merging experimental results}, top rows).
Then, it will tend to leave space for the red vehicle to merge into the right lane. If the red vehicle is aware of this courteous behavior, it merges into the right lane accordingly. However, if it wrongfully anticipates stubborn behavior, it will not attempt to merge into the right lane and keep driving next to the other two vehicles. When approaching the white vehicle, it settles to the strategy of decelerating and merging behind the yellow vehicle. This is suboptimal but not fatal, and accordingly marked in orange.

On the other hand, consider the scenario in which the yellow vehicle is stubborn (\cref{fig:two-dimensional GPG merging experimental results}, bottom rows). If the red vehicle is aware of this behavior, it does not attempt to merge into the right lane until the yellow vehicle has passed. However, under the \textit{courteous} belief, it plans a future open-loop trajectory which merges into the other lane, (wrongfully) assuming that the yellow vehicle will make space for this manoeuvre. During this manoeuvre, collision avoidance constraints are violated, leading to undesired and dangerous behavior (marked in red).

It is apparent that wrong parameter estimates may prohibit safe deployment of a control scheme of this type. 
We now perform the same simulations using the online learning scheme \eqref{mini:human_parameter_estimation}, which continuously updates the estimate of $\params[\yellowplayer]$, with $L = 1$ and regularization $\xi = 0.5$. At the first time step, as there are no observations available yet, the open-loop sequence predicted by the red vehicle is equal to the open-loop sequence obtained previously.
On the other hand, the behavior at the 40\textsuperscript{th} time step has significantly improved and the previously observed suboptimal/dangerous maneuvers no longer occur.
\captionsetup[subfigure]{labelformat=empty}
\begin{figure}[ht]
	\centering
	\begin{tabular}{wl{0.20\columnwidth} >{\centering\arraybackslash}m{0.69\columnwidth}}
		\toprule
		\hspace{-0.3cm} \scriptsize Belief & \hspace{-2.2cm} \normalsize Courteous vehicle\\
		\midrule
	\end{tabular}
	\begin{subfigure}[t]{0.49\columnwidth}
		\makebox[0.06\columnwidth][r]{\makebox[10pt]{\raisebox{15pt}{\rotatebox[origin=c]{90}{\scriptsize Stubborn}}}}%
		\adjincludegraphics[width = \columnwidth-1pt, trim={0 0 {.5\width} 0},clip,cfbox=green 1pt 0pt]{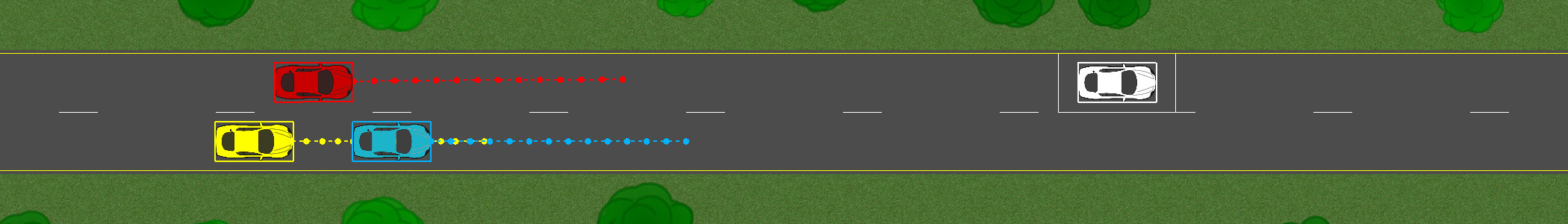}
	\end{subfigure}
	\hspace{-0.015\columnwidth}
	\begin{subfigure}[t]{0.49\columnwidth}
		\adjincludegraphics[width = \columnwidth-1pt, trim={0 0 {.5\width} 0},clip, cfbox=green 1pt 0pt]{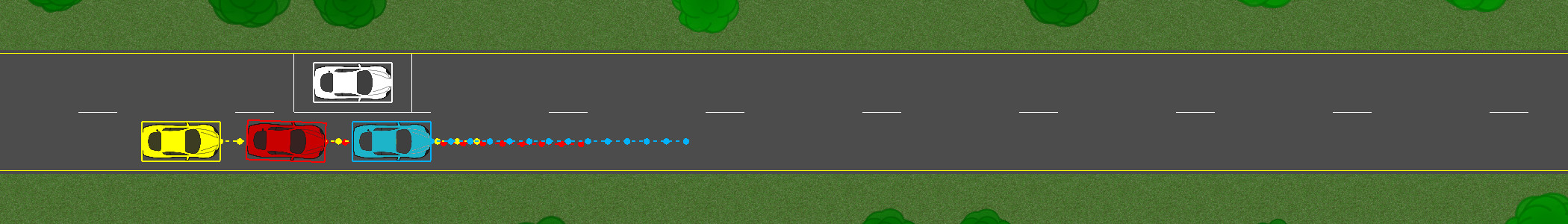}
	\end{subfigure}
	\vspace{-0.3cm}
	\phantom{\tiny abc}
	\begin{tabular}[ht]{m{0.20\columnwidth} >{\centering\arraybackslash}m{0.69\columnwidth}}
		\toprule
		\hspace{-0.3cm} \scriptsize Belief & \hspace{-2.2cm} \normalsize Stubborn vehicle\\
		\midrule
	\end{tabular}
	\begin{subfigure}[t]{0.49\columnwidth}
		\makebox[0.06\columnwidth][r]{\makebox[10pt]{\raisebox{15pt}{\rotatebox[origin=c]{90}{\scriptsize Courteous}}}}%
		\adjincludegraphics[width = \columnwidth-1pt, trim={0 0 {.5\width} 0},clip, cfbox=green 1pt 0pt]{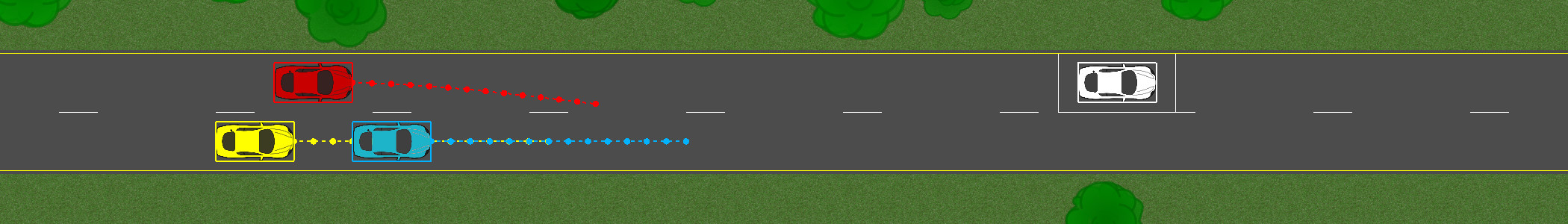}
		\caption{Time step 1}
	\end{subfigure}
	\hspace{-0.015\columnwidth}
	\begin{subfigure}[t]{0.49\columnwidth}
		\adjincludegraphics[width = \columnwidth-1pt, trim={0 0 {.5\width} 0},clip, cfbox=green 1pt 0pt]{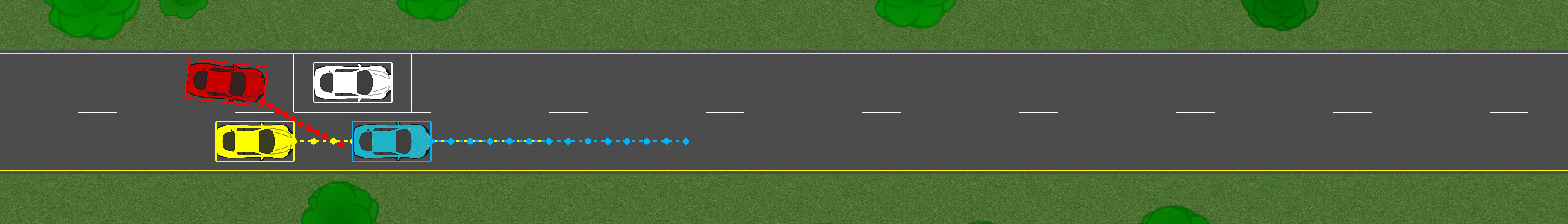}
		\caption{Time step 40}
	\end{subfigure}
	\caption{The planned trajectories at respectively the 1\textsuperscript{st} and 40\textsuperscript{th} time step of the game when using online learning.}
	\vspace{-0.25cm}
	\label{fig:two-dimensional GPG_merging experimental results belief}
\end{figure}

\ifArxiv
\cref{fig:merging parameter estimates} shows the evolution of the parameter estimate $\params{\yellowplayer}_1$ over time.
\fi
As discussed in \Cref{sec:online_learning}, we do not observe convergence to the true parameter value (which is to be expected whenever the control signals are not sufficiently exciting). Nevertheless, parameter estimates are obtained which better explain previous observations, resulting in significantly improved behavior.
\ifArxiv
\begin{figure}
	\centering
	\begin{tikzpicture}
	\begin{groupplot}[
	group style={group size=2 by 1},
	width=\columnwidth/2, 
	height=3cm,
	]
	\nextgroupplot[title=Courteous vehicle, xlabel=Time step, ylabel=$\params{\yellowplayer}_1$,  ytick={0, 5, 10}]
	\addplot[thick] file 
	{experiments/merging_scenario_ECC/True-True-True/iter/belief_0_1_0.dat};
	
	\addplot[thick, dashed] file 
	{experiments/merging_scenario_ECC/True-True-False/iter/belief_0_1_0.dat};
	
	\nextgroupplot[title=Stubborn vehicle, xlabel=Time step, ytick={0, 5, 10}]
	\addplot[thick] file 
	{experiments/merging_scenario_ECC/True-False-False/iter/belief_0_1_0.dat};

	\addplot[thick, dashed] file 
	{experiments/merging_scenario_ECC/True-False-True/iter/belief_0_1_0.dat};
	\end{groupplot}
	\end{tikzpicture}
	\caption[The estimate of the human parameter.]{The estimate of $\params{\yellowplayer}_1$ over time, where the full lines denote experiments with correct initial estimate and the dashed lines denote experiments with incorrect initial estimate.}
	\label{fig:merging parameter estimates}
\end{figure}
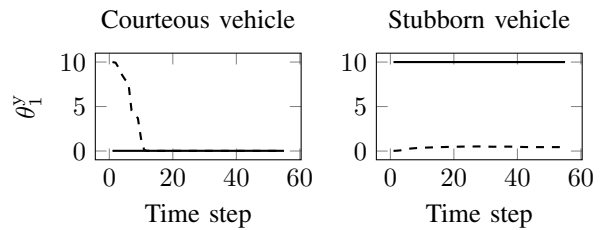
\fi 
This conclusion is corroborated by \Cref{tab:two-dimensional GPG merging potential and violation}, which provides the numerical values for the closed-loop potential and the maximum constraint violation. 
For both experiments with incorrect initial belief, introducing the parameter estimation scheme results in significantly improved closed-loop potential, approaching the value for the respective experiment with correct initial belief. Furthermore, the maximum constraint violation is smaller than the imposed tolerance $10^{-2}$. Note that for both experiments with correct initial belief, the closed-loop potential and maximum constraint violation are (approximately) unaffected by the learning methodology.
Finally, \Cref{tab:two-dimensional GPG merging computation times} provides the computation times for the controlled vehicle, illustrating the real-time capabilities of the overall scheme, as the maximum computation time for solving the GPG and performing online learning combined are below our sampling period $\Ts = \SI{0.2}{\second}$.
\begin{table}
	\vspace{0.15cm}
	\caption{The closed-loop potential and maximum constraint violation.}
	\label{tab:two-dimensional GPG merging potential and violation}
	\centering
	\begin{filecontents*}{experiments/merging_scenario_ECC/potential_violation.csv}
Behavior,Belief,Simple,Learning,Simple,Learning
Courteous,Courteous,   36.1757919623,   36.1734551045,    0.0000000000,    0.0000000000
Courteous,Stubborn,   74.5530504414,   35.9573532389,    0.0000000000,    0.0000000000
Stubborn,Courteous,  125.4973078423,   44.4726490458,    3.6669081268,    0.0086314009
Stubborn,Stubborn,   44.1246099017,   44.1251950346,    0.0095409141,    0.0097823581
\end{filecontents*}
\pgfplotstabletypeset[
	col sep = comma,
	fixed, fixed zerofill, precision = 3,
	string replace*={_}{\textsubscript},
	every head row/.style={
		before row={
			\toprule
			Behavior & Belief & \multicolumn{2}{c}{Potential} & \multicolumn{2}{c}{Violation}\\
			\cmidrule(lr){3-4}\cmidrule(lr){5-6}
		},
		after row=\midrule,
	},
	every last row/.style={after row=\bottomrule},
	display columns/0/.style={string type,column name={}, column type = {l},
		multicolumn names},
	display columns/1/.style={string type,column name={}, column type = {l},
		multicolumn names},
	display columns/2/.style={fixed, precision = 2, 1000 sep={ }},
	display columns/3/.style={fixed, precision = 2, 1000 sep={ }},
	display columns/4/.style={fixed, precision = 3, 1000 sep={ }, column name={Simple}},
	display columns/5/.style={fixed, precision = 3, 1000 sep={ }, column name={Learning}},
	]{experiments/merging_scenario_ECC/potential_violation.csv}
\end{table}
\begin{table}
	\caption{The computation times of the red vehicle.}
	\label{tab:two-dimensional GPG merging computation times}
	\centering
	\begin{filecontents*}{experiments/merging_scenario_ECC/timing_learning_reduced.csv}
Behavior,Belief,Max,Mean,Max,Mean
Courteous,Courteous,    0.0129430000,    0.0025609818,    0.0055210000,    0.0017259815
Courteous,Stubborn,    0.0935390000,    0.0127627455,    0.1030110000,    0.0255024815
Stubborn,Courteous,    0.0367760000,    0.0117785636,    0.1023450000,    0.0427574259
Stubborn,Stubborn,    0.0599850000,    0.0130962182,    0.1019030000,    0.0306740185
\end{filecontents*}
\pgfplotstabletypeset[
	col sep = comma,
	fixed, fixed zerofill, precision = 3,
	string replace*={_}{\textsubscript},
	every head row/.style={
		before row={
			\toprule
			Behavior & Belief & \multicolumn{2}{c}{GPG [s]} & \multicolumn{2}{c}{Update $\params{\yellowplayer}$ [s]}\\
			\cmidrule(lr){3-4}\cmidrule(lr){5-6}
		},
		after row=\midrule,
	},
	every last row/.style={after row=\bottomrule},
	display columns/0/.style={string type,column name={}, column type = {l},
		multicolumn names},
	display columns/1/.style={string type,column name={}, column type = {l},
		multicolumn names},
	display columns/2/.style={column name={Max}},
	display columns/3/.style={column name={Avg}},
	display columns/4/.style={column name={Max}},
	display columns/5/.style={column name={Avg}},
	]{experiments/merging_scenario_ECC/timing_learning_reduced.csv}
\end{table}

\section{Conclusion}
We presented an interaction-aware MPC strategy for 
automated driving under shared collision avoidance 
constraints by formulating the problem as a \ac{GPG}, which depends on implicit preferences and constraints 
of surrounding drivers. We present a simple but effective 
scheme which allows the controller to learn these preferences online using observed behavior. 
Through numerical simulations, we have illustrated the benefits of the learning scheme, as well as the potential real-time capabilities of the proposed methodology. 

\scriptsize
\bibliography{zotero_library}
\bibliographystyle{ieeetr}
\end{document}